\documentclass[a4paper]{article}

\usepackage{delarray,verbatim,enumerate,a4}
\usepackage{amsmath,amsthm,amstext,amsbsy,amssymb,amsfonts,amscd}
\usepackage[frenchb]{babel}
\usepackage{aeguill}

\newtheorem{Th}{Th\'eor\`eme}[]

\def\PreuveTh{\smallskip\noindent {\it Preuve du Théorème.~}}

\font\teneufm=eufm10
\font\seveneufm=eufm7
\font\fiveeufm=eufm5
\newfam\gothfam
\textfont\gothfam=\teneufm \scriptfont\gothfam=\seveneufm
\scriptscriptfont\gothfam=\fiveeufm
\def\goth{\fam\gothfam}

							\def\QQ{\mathbb Q}	
\def\NN{\mathbb N}			\def\ZZ{\mathbb Z}
\def\F2{\mathbb{F}_2}		\def\Z2{\mathbb{Z}_2}	
\def\Zl{\mathbb{Z}_\ell} 	\def\Ql{\mathbb{Q}_\ell}

 				\def\U{\mathcal  U}
  		\def\C{\mathcal  C}		
 	  	\def\Cl{\mathcal  C \ell}
\def\E{\mathcal  E}		

				\def\l{{\goth l}}

\def\Gal{\operatorname{Gal}}

\begin{document}

\title{\LARGE\bf Note sur la conjecture de Leopoldt}

\author{ Jean-François {\sc Jaulent} }
\date{}
\maketitle

{\small
\noindent{\bf Résumé.} Nous montrons qu'un corps de nombres de degré arbitraire satisfait la conjecture de Leopoldt dès lors qu'il est peu ramifié en un sens bien précis.
}

\

{\small
\noindent{\bf Abstract.} We prove that number fields with arbitrary  degree but weak ramification satisfy the Leopoldt conjecture on the $\ell$-adic rank of the group of units. 
}
\bigskip\bigskip

\bigskip

\noindent{\large \bf Introduction}

\medskip

La conjecture de Leopoldt  postule que pour chaque nombre premier $\ell$, le rang $\ell$-adique du groupe des unités globales $E_K$  d'un corps de nombres est égal au $\ZZ$-rang de ce groupe, c'est à dire, conformément au théorème de Dirichlet, à la somme $r_K+c_K-1$ des nombres de places réelles ou complexes de ce corps diminuée de 1.
\smallskip

Diverses approches ont été tentées, qui la prouvent dans certains cas~:\smallskip

Les {\em méthodes algorithmiques} permettent naturellement de constater qu'un corps de nombres donné $K$ satisfait la conjecture de Leopoldt, dès lors que l'on sait {\em effectivement} faire des calculs dans $K$, i.e. déterminer pratiquement un système d'unités fondamentales de l'anneau des entiers. Il suffit alors, en effet,  de calculer un certain régulateur $\ell$-adique avec une précision convenable, pour établir qu'il n'est pas nul. Ces méthodes permettent donc  de vérifier la conjecture pour tout corps $K$ {\em donné} et un premier $\ell$ {\em fixé}, dès lors que le degré $n_K$ de $K$ n'est pas trop grand.\smallskip

Les {\em méthodes algébriques} donnent des conditions suffisantes (mais non nécessaires) à sa validité. La théorie d'Iwasawa montre ainsi que le défaut de la conjecture de Leopoldt pour un premier fixé $\ell$ dans la $\Zl$-extension cyclotomique $K^c=\cup_{n\in\NN} K_n$ d'un corps de nombres $K$ est majoré par l'invariant lambda attaché à la limite projective $\C_{K^c}=\varprojlim \Cl'_{K_n}$ des $\ell$-groupes de $\ell$-classes des étages finis de la tour (cf. e.g. \cite{W}). Sous l'hypothèse  de trivialité $\C_{K^c}=1$, la conjecture de Leopoldt est donc vérifiée à tous les étages de la tour. Cette observation permet ainsi de construire des familles infinies de corps de nombres qui la satisfont (cf. e.g. \cite{G, J2, J3}).

\smallskip

Les {\em méthodes transcendantes} s'appuient principalement sur les résultats d'indé\-pendance de logarithmes de nombres algébriques et utilisent très peu les propriétés arithmétiques des unités. C'est pourquoi elles se généralisent sans peine à n'importe quel sous-groupe de type fini du groupe multiplicatif $K^\times$ (cf. \cite{J1}). Appliqué dans un contexte galoisien, le théorème de Baker-Brumer permet ainsi d'établir la conjecture de Leopoldt dans le cas où le corps considéré $K$ est une extension abélienne de $\QQ$ ou d'un corps quadratique imaginaire $k$ et dans quelques autres situations (en particulier dès que l'algèbre de Galois $\Ql [\Gal (K/k)]$ est un produit direct de corps).

\bigskip
\noindent{\large \bf Théorème principal}
\medskip

Le but de la présente note est de présenter une quatrième approche basée sur les résultats analytiques d'Odlysko, Poitou et Serre sur les minorations de discriminant. Son résultat principal affirme qu'un corps de nombres de degré arbitraire satisfait la conjecture de Leopoldt pour tous les premiers $\ell$ dès lors qu'il n'est pas trop ramifié. \par
Il s'énonce comme suit~:

\begin{Th}
Soit $F$ un corps de nombres de degré $n_F$ et de discriminant absolu $d_F$. Sous la condition

\centerline{$|d_F|^{1/n_F} \le 22,3$,}\smallskip

\noindent le corps $F$ vérifie la conjecture de Leopoldt pour tous les nombres premiers $\ell$.\par
 Et la même conclusion vaut encore, si l'on admet la conjecture de Riemann généralisée, sous la condition plus faible~:\smallskip

\centerline{$|d_F|^{1/n_F} \le 44,7$.}
\end{Th}

\PreuveTh Elle est très simple~: nous allons procéder par contraposée en supposant que le corps considéré $F$ ne satisfait pas la conjecture de Leopoldt pour un premier donné $\ell$~; autrement dit que qu'il existe un élément non trivial $\varepsilon$ dans le tensorisé  $\ell$-adique $\E_F=\Zl\otimes_\ZZ E_F$ du groupe des unités globales de $F$ d'image locale triviale dans le groupe $\prod_{\l | \ell} \U_\l$ des unités semi-locales attaché aux places $\l$ de $F$ au-dessus de $\ell$ (cf. \cite{J3}).
\par
 Pour chaque entier  $n \ge 1$, écrivons $\varepsilon = \varepsilon_n \,\eta_n^{\ell^n}$ avec  $\eta_n$ dans $\E_F$ et $\varepsilon_n$ dans $E_F$ suffisamment proche de 1 (disons de valuation $\nu_\l(\varepsilon_n-1)>\ell^n+\kappa$)  en chacune des places $\l$ divisant $\ell$, de sorte que nous puissions définir sa racine $\ell^n$-ième $\xi_n$ par~:\smallskip

\centerline{$\xi_n = \big(1+(\varepsilon-1)\big)^{1/\ell^n} = \sum_{k=0}^\infty \binom{1/\ell^n}{k}(\varepsilon-1)^k$.}\smallskip

Considérons alors pour chaque $n>0$ l'extension $F_n$ de $F$ engendrée par $\xi_n$.
\par 
Observons que le degré $[F_n:F]$ tend vers l'infini avec $n$ ; que l'extension $F_n/F$ est non ramifiée en dehors de $\ell$, puisque $\varepsilon_n$ est une unité ; et qu'elle est par construction complétement décomposée aux places au-dessus de $\ell$. En particulier que $F_n/F$ est non ramifiée à toutes les places finies et que la quantité $|d_{F_n}|^{1/[F_n:\QQ]}$, construite sur le discriminant $d_{F_n}$ de $F_n$,  est indépendante de $n$ et égale à $|d_{F}|^{1/[F:\QQ]} =|d_{F}|^{1/n_F}$.\par
Cela étant, les minorations asymptotiques obtenues via des considérations analytiques par Odlysko, Poitou et Serre nous donnent l'inégalité~:
\smallskip

\centerline{$|d_F|^{1/n_F} > 22,3$~;}\smallskip

\noindent et même l'inégalité plus forte~:\smallskip

\centerline{$|d_F|^{1/n_F} > 44,7$,}\smallskip

\noindent sous l'hypothèse de Riemann généralisée (cf. e.g. \cite{M}). D'où le résultat annoncé.


\def\refname{\small{\sc  Références}}

\bigskip\noindent
{\small
\begin{tabular}{l}
{Jean-Fran\c cois {\sc Jaulent}}\\
Institut de Math{\'e}matiques de Bordeaux \\
Université {\sc Bordeaux 1} \\
351, cours de la lib{\'e}ration\\
F-33405 {\sc Talence} Cedex\\
courriel : Jean-Francois.Jaulent@math.u-bordeaux1.fr 
\end{tabular}
}


{\small

\begin{thebibliography}{tttt}



  
\bibitem[Gr]{G} {\sc G. Gras,}
{\it Class Field Theory},
Springer Monographs in Mathematics (2003).

\bibitem[Ja$_1$]{J1} {\sc J.-F. Jaulent},
{\it Sur l'indépendance  $\ell$-adique de nombres algébriques}, 
J. Numb. Th. {\bf 20} (1985), 149--158 .
 
\bibitem[Ja$_2$]{J2} {\sc J.-F. Jaulent},
{\it Sur les conjectures de Leopoldt et de Gross}, 
Actes des Journées Arithmétiques de Besançon, Astérisque 
{\bf 147-148} (1987), 107--120.

\bibitem[Ja$_3$]{J3}
{\sc J.-F. Jaulent},
{\em Théorie $\ell$-adique globale du corps de classes},
 J. Théor.  Nombres Bordeaux {\bf 10} (1998),   355--397.
 
 \bibitem[Ma]{M}
{\sc J. Martinet},
{\em Tours de corps de classes et estimations de discriminants},
 Invent. Math. {\bf 44} (1978),   65--73.
 
 \bibitem[Wa]{W} {\sc L. Washington,}
{\it Introduction to cyclotomic fields},
Graduate Texts in Mathematics {\bf 83}. Springer-Verlag, New York, 1997. xiv+487 pp


\end{thebibliography}
}
 \end{document}